# STABILITY AND THE LYAPOUNOV EXPONENT OF THRESHOLD AR-ARCH MODELS

By Daren B. H. Cline and Huay-min H. Pu

*Texas A&M University*

The Lyapounov exponent and sharp conditions for geometric ergodicity are determined of a time series model with both a threshold autoregression term and threshold autoregressive conditional heteroscedastic (ARCH) errors. The conditions require studying or simulating the behavior of a bounded, ergodic Markov chain. The method of proof is based on a new approach, called the *piggyback method*, that exploits the relationship between the time series and the bounded chain.

The piggyback method also provides a means for evaluating the Lyapounov exponent by simulation and provides a new perspective on moments, illuminating recent results for the distribution tails of GARCH models.

**1. Introduction.** Modeling the stochastic volatility of econometric and other time series with autoregressive conditional heteroscedastic (ARCH and GARCH) type models has proved to be very successful. This effort recently has been extended to include additional nonlinearity such as threshold (G)ARCH models introduced by Glosten, Jagannathan and Runkle (1993), Rabemanajara and Zakoian (1993) and Zakoian (1994), as well as the addition of autoregression components, which may also be nonlinear [Li and Li (1996), Liu, Li and Li (1997), Lu (1998), Hwang and Woo (2001), Lu and Jiang (2001) and Lanne and Saikkonen (2004)]. For nonlinear models in particular, it is quite usual to make very strong assumptions about the parameters and/or nonlinearity in order to ensure existence of a stationary model. To avoid such strong assumptions, this paper introduces an approach that will give definitive results about stationarity and can offer insight to the moments of the stationary distribution. The method is applied here to threshold extensions of ARCH models with an autoregressive component.

The threshold autoregressive ARCH (TAR-ARCH) model can be expressed as

$$\xi_t = a(\xi_{t-1}, \ldots, \xi_{t-p}) + b(\xi_{t-1}, \ldots, \xi_{t-p})e_t, \tag{1.1}$$

where $a$ and $b$ are finitely piecewise continuous functions and $\{e_t\}$ is an i.i.d. error sequence. The state vector for the time series is $X_t = (\xi_t, \ldots, \xi_{t-p+1})$. The nonlinear autoregression function $a(x)$ is continuous on individual, connected subregions of $\mathbb{R}^p$; the boundaries of these regions are called thresholds, hence the nomenclature for the model. Frequently, $a(x)$ is assumed to be linear on each of these regions. Likewise, this model has a state dependent conditional variance, $b^2(x) = \text{var}(\xi_t | X_{t-1} = x)$ if $\text{var}(e_t) = 1$, which typically is of the order of magnitude of $\|x\|^2$. This provides the conditional heteroscedasticity (ARCH) behavior and for our purposes it is also assumed to have thresholds.

Since the time series (1.1) is embedded into $\{X_t\}$, which is a Markov chain, it will have a stationary distribution when the Markov chain is ergodic. If $a$ and $b$ are sufficiently smooth (e.g., Lipschitz continuous) this may be verified using dynamical systems arguments [Chan and Tong (1985), Letac (1986), Tong (1990) and Diaconis and Freedman (1999)]. Here, however, the discontinuities of $a$ and $b$, as well as their nonlinearity, complicate every aspect of the argument. Another approach, which works well for nonthreshold ARCH models, involves reexpressing the model as a random coefficients or stochastic difference equation model [Brandt (1986), Engle and Bollerslev (1986), Bougerol (1987), Bougerol and Picard (1992a, b) and many others]. This approach can even provide stationary representations and tail behaviors of the stationary distributions [Basrak, Davis and Mikosch (2002)], but it will not work here since the "random coefficients" would not be independent of the past values of the process. Indeed the approach assumes that the random coefficients are already known to be a stationary process. Our approach will overcome these obstacles.

Both the dynamical systems approach and the random coefficients approach are intricately tied to the concept of a Lyapounov exponent. This is a notion readily apparent in such systems but largely bypassed in the general literature on ergodicity of Markov processes and the like. An appropriate definition for the Lyapounov exponent in our context is

$$\liminf_{n \to \infty} \limsup_{\|x\| \to \infty} \frac{1}{n} E\left(\log\left(\frac{\|X_n\|}{\|X_0\|}\right) \Big| X_0 = x\right).$$

This constant measures the "geometric drift" of the process when it is large (and less subject to local perturbations). If this constant is negative, then the process is stable in the sense that a drift condition for ergodicity is readily identifiable [cf. Meyn and Tweedie (1993), Theorem 15.3.7]. Our approach in this effort will both identify the Lyapounov exponent and show that it is the critical value for determining stability.

As is standard, we assume the following on the error term, ensuring $\{X_t\}$ is an aperiodic, Lebesgue-irreducible $T$-chain on $\mathbb{R}^p$ [cf. Meyn and Tweedie (1993) and Cline and Pu (1998)].

ASSUMPTION A.1. The distribution of $e_t$ has Lebesgue density $f$ which is locally bounded away from 0. Also, $b$ is positive, locally bounded and locally bounded away from 0.

We also have the following moment assumption.

ASSUMPTION A.2. $\sup_u (1 + |u|) f(u) < \infty$ and $E(|e_1|^{r_0}) < \infty$ for some $r_0 > 0$.

In this paper we will employ further constraints on the functions $a$ and $b$. These conform, however, to all standard uses of threshold autoregression and ARCH models. The first such condition is next. Except when indicated otherwise, $\|\cdot\|$ is the Euclidean norm though this is not strictly necessary.

ASSUMPTION A.3. $a(x)/(1+\|x\|)$ and $b(x)/(1+\|x\|)$ are bounded.

The model thus has the functional coefficient AR-ARCH (FCAR-ARCH) representation

$$
\begin{aligned}
\xi_t = {} & a_0(X_{t-1}) + \sum_{i=1}^{p} a_j(X_{t-1}) \xi_{t-i} \\
& + \left( b_0^2(X_{t-1}) + \sum_{i=1}^{p} b_j^2(X_{t-1}) \xi_{t-i}^2 \right)^{1/2} e_t,
\end{aligned}
\tag{1.2}
$$

where $a_0, \ldots, a_p, b_0, \ldots, b_p$ are bounded functions on $\mathbb{R}^p$ [e.g., $a_0(x) = \frac{a(x)}{1+|x_1|+\cdots+|x_p|}$ and $a_i(x) = \frac{\operatorname{sgn}(x_i) a(x)}{1+|x_1|+\cdots+|x_p|}$, $i = 1, \ldots, p$]. Although the FCAR-ARCH representation is not unique, it seems to be a good starting point for (nonparametric) modeling. Stability conditions based on such a representation may not be sharp, however, and this is one motivation for the present work.

Standard threshold models consider $a$ to be piecewise linear and $b^2$ to be piecewise quadratic. More generally, it is not unreasonable to assume that $a(x)$ and $b(x)$ are asymptotically homogeneous. We state this assumption next. Let $\Theta = \{\theta \in \mathbb{R}^p : \|\theta\| = 1\}$.

ASSUMPTION A.4. There exist $x_* \in \mathbb{R}^p$ and bounded functions $a_*$ and $b_*$ on $\Theta$ such that

$$\lim_{w\to\infty} \sup_{\theta\in\Theta} \left| \frac{a(x_* + w\theta)}{w} - a_*(\theta) \right| = 0 \quad \text{and} \quad \lim_{w\to\infty} \sup_{\theta\in\Theta} \left| \frac{b(x_* + w\theta)}{w} - b_*(\theta) \right| = 0.$$

This actually implies Assumption A.3, but it is convenient to state both. It follows easily from Assumption A.4 that $a$ and $b$ may be decomposed as $a(x) = a_0(x) + a^*(x - x_*)$ and $b(x) = b_0(x) + b^*(x - x_*)$, $x = (x_1, \ldots, x_p) \in \mathbb{R}^p$, where $a^*$ and $b^*$ are homogeneous [in fact, $a^*(x) = a_*(x/\|x\|)\|x\|$ and $b^*(x) = b_*(x/\|x\|)\|x\|$] and $a_0$ and $b_0$ are locally bounded functions such that $a_0(x) = o(\|x\|)$ and $b_0(x) = o(\|x\|)$ as $\|x\| \to \infty$. For simplicity, we henceforth assume $x_* = 0$.

By removing the nonhomogeneous terms from the model, we define a homogeneous version of the time series,

$$(1.3) \qquad \xi_t^* = a^*(X_{t-1}^*) + b^*(X_{t-1}^*)e_t,$$

where $X_t^* = (\xi_t^*, \ldots, \xi_{t-p+1}^*)$. We intend to show that, under simple conditions, the stability of $\{X_t\}$ is related to that of $\{X_t^*\}$. Furthermore, letting $\theta_t^* = X_t^*/\|X_t^*\|$, we can describe the stability conditions in terms of the collapsed process $\{\theta_t^*\}$. Note that $\{X_t^*\}$ is a Markov chain on $\mathbb{R}_0^p = \mathbb{R}^p \setminus \{0\}$ and, due to the homogeneity of $a^*$ and $b^*$, $\{\theta_t^*\}$ is a Markov chain on $\Theta = \{\theta : \|\theta\| = 1\}$.

When $p = 1$, the latter chain obviously is a two-state chain and we examine this special case explicitly in Example 4.1. Much more interesting is the case $p > 1$.

Even though we are assuming $b$ to be locally bounded away from 0, $b^*$ need not be. For example, for the standard ARCH model $b^*(x) = (b_1^2 x_1^2 + \cdots + b_p^2 x_p^2)^{1/2}$, where some $b_i$ may be 0. This requires an additional assumption and some care. We define $B_0$ to be the set on which $\liminf_{y \to x} b^*(y) = 0$ and $H_0 = \{x : \min_i |x_i| = 0\}$.

ASSUMPTION A.5. If $p > 1$, $\max(|a^*(x)|, b^*(x))$ is locally bounded away from 0 on $\mathbb{R}_0^p$ and $b^*$ is locally bounded away from 0 on $\mathbb{R}_\#^p = \mathbb{R}^p \setminus H_0$ (hence $B_0 \subset H_0$).

Our final assumption refers to the threshold-like nature of the model, that is, the piecewise continuity of $a^*$ and $b^*$.

ASSUMPTION A.6. If $p > 1$, there exist $(p-1)$-dimensional homogeneous hyperplanes $H_1, \ldots, H_m$ such that $a^*$ and $b^*$ are continuous off $\bigcup_{j=1}^m H_j$.

Assumption A.6 implies that the character of the model depends on the signs of linear combinations of the components of $X_t$. This is similar to and includes the delay assumption most authors employ.

Define $z(x,u) = a^*(x) + b^*(x)u$ for $x \in \mathbb{R}^p$ and $u \in \mathbb{R}$, and $w(x,u) = \|(z(x,u), x_1, \ldots, x_{p-1})\|$. Then $\xi_t^* = z(X_{t-1}^*, e_t)$ and $\theta_t^* = (z(\theta_{t-1}^*, e_t), \theta_{t-1,1}^*, \ldots, \theta_{t-1,p-1}^*)/w(\theta_{t-1}^*, e_t)$.

Henceforth, let $\mu_p$ be Lebesgue measure on $\mathbb{R}^p$ and let $\mu_\Theta$ be the surface area measure on $\Theta$. That is, $\mu_\Theta(A) = p\mu_p(\{x : x/\|x\| \in A, \|x\| \leq 1\})$.

Section 2 discusses a direct approach to verifying stability that has its roots in the literature for both parametric and nonparametric models. We present it because of its simplicity and to contrast it with our results in Section 3. The main results, which are in Section 3, provide a sharp condition for stability and identify the Lyapounov exponent for $\{X_t\}$. The approach is based on a method of proof we have named the *piggyback method*. Section 4 examines four cases explicitly. The final section contains the proofs.

**2. A simple sufficient condition for FCAR-ARCH models.** In the spirit of numerous other authors, we present a simple condition for geometric ergodicity applicable to general functional coefficient models. This condition is based on bounding the FCAR-ARCH representation and therefore is simple to apply, but it ultimately proves to be too rudimentary to establish the parameter space definitively for more specific models. The approach has roots in conditions for FCAR models, starting with Chan and Tong (1986), and culminating recently in conditions for models with both a nonlinear autoregression and an ARCH or bilinear component [Lu (1998), Lu and Jiang (2001) and Ferrante, Fonseca and Vidoni (2003)]. Our theorem and corollary below encompass these latest results.

First, however, we briefly recall the Foster–Lyapounov drift condition for establishing geometric ergodicity of a Markov chain as developed by Meyn and Tweedie [(1993), Theorem. 16.0.1], namely the chain is aperiodic, $\phi$-irreducible, $T$-continuous and there exist $K < \infty$, $\beta_0 < 1$ and a test function $V(x) \geq 1$ such that $V(x) \to \infty$ as $\|x\| \to \infty$ and

$$E(V(X_1)|X_0 = x) \leq K + \beta_0 V(x) \qquad \text{for all } x.$$

This condition in fact establishes that $\{X_t\}$ is $V$-uniformly ergodic: with $\Gamma$ as the stationary distribution and $\beta_1 \in (\beta_0, 1)$,

$$(2.1) \qquad \sup_A |P(X_n \in A|X_0 = x) - \Gamma(A)| \leq M\beta_1^n V(x) \qquad \text{all } n, x,$$

for some finite $M$. The further advantage of knowing something about the function $V$ is twofold. First, $\int V(x)\Gamma(dx) < \infty$ [cf. Meyn and Tweedie (1993), Theorem 14.0.1], implying that if $|g(x)| \leq V(x)$, then $\int |g(x)|\Gamma(dx) < \infty$ and $\frac{1}{n}\sum_{t=1}^n g(X_t) \to \int g(x)\Gamma(dx)$ almost surely. Second, a central limit theorem holds for $\sum_{t=1}^n g(X_t)$ if $|g(x)| \leq V(x)^{1/2}$ [cf. Meyn and Tweedie (1993), Theorem 17.0.1]. In particular, an appropriate test function will ensure the existence of moments.

THEOREM 2.1. *Assume Assumptions* A.1–A.3. *If there exist nonnegative $c_i$ with $\sum_{i=1}^p c_i < 1$, $r > 0$ and $K < \infty$ such that*

$$E(|\xi_1|^r | X_0 = x) \leq K + \sum_{i=1}^p c_i |x_i|^r, \tag{2.2}$$

*for all $x = (x_1, \ldots, x_p) \in \mathbb{R}^p$, then $\{X_t\}$ is $V$-uniformly ergodic with $V(x) = 1 + \sum_{i=1}^p d_i |x_i|^r$ for some positive $d_1, \ldots, d_p$. Furthermore, the stationary distribution has finite $r$th moment.*

The theorem can be applied when one is given a specific representation of the model. The following corollary does so, extending Lu (1998) and Lu and Jiang (2001), who handle the cases $r = 2$ and $r = 1$, respectively.

COROLLARY 2.2. *Assume Assumptions* A.1–A.3 *and that there are nonnegative constants $a_0, \ldots, a_p$, $b_0, \ldots, b_p$ such that*

$$|a(x)| \leq a_0 + a_1|x_1| + \cdots + a_p|x_p| \quad \text{and} \quad b(x) \leq (b_0^2 + b_1^2 x_1^2 + \cdots + b_p^2 x_p^2)^{1/2}.$$

*Each of the following implies (2.2) holds with $\sum_{i=1}^p c_i < 1$:*

(i) $r \leq 1$ *and* $\sum_{i=1}^p (a_i^r + b_i^r E(|e_t|^r)) < 1$.

(ii) *Either* $1 < r < 2$ *and the error distribution is symmetric about* 0, *or $r = 2$ and $E(e_1) = 0$, and* $(\sum_{i=1}^p a_i)^r + \sum_{i=1}^p b_i^r E(|e_t|^r) < 1$.

Although no particular FCAR-ARCH representation is assumed, it seems one cannot be avoided in checking the condition. Its advantage is in its simplicity. The conditions presented in the next section, on the other hand, are sharp.

**3. The piggyback method and the main results.** The *piggyback method* is the name we have given to a new approach for determining stability of nonlinear time series. This approach has proven useful for some rather simple models that do not behave the same as their deterministic "skeleton." See Example 3.2 of Cline and Pu (1999) and Example 4.2 of Cline and Pu (2002). See also Gourieroux and Monfort (1992) for a simple example. We anticipate that this approach will prove to have its greatest use for models such as those studied here, namely for models in which the error terms are state dependent and thereby have a profound effect on stability. The name "piggyback" derives from the use of a Foster–Lyapounov drift condition that piggybacks on the stability of another Markov process. This second process is related to a process embedded in the first, is more basic and has a much more obvious criterion for stability.

For the threshold AR-ARCH model we are investigating, the related process is $\{\theta_t^*\}$ (the *collapsed* process) which, being a process on compact $\Theta$, is clearly bounded. With minimal conditions for irreducibility and so on, therefore, it is uniformly ergodic. We note that $w(\theta_0^*, e_1)$ (see Section 1) is a surrogate for the one-step change in magnitude of the AR-ARCH model, namely $\|X_1\|/\|X_0\|$ when $\|X_0\|$ is large. Intuitively then, the Markov chain $\{X_t\}$ should be stable if the expected value of $\log(w(\theta_0^*, e_1))$, taken relative to the stationary distribution of $\{\theta_t^*\}$, is less than 0. The trick is to construct a function to verify this, especially in view of the discontinuous coefficient functions, and this is where the piggyback method is both extremely helpful and somewhat intricate.

We now state the stability theorems for these two processes. An outline of the method follows the statement of Theorem 3.3.

THEOREM 3.1. *Assume Assumptions* A.1–A.5. *Then $\{\theta_t^*\}$ is a $\mu_\Theta$-irreducible, uniformly ergodic T-chain on $\Theta$ with stationary distribution $\Pi$ satisfying*

$$(3.1) \qquad \int_\Theta \int_\mathbb{R} |z(\theta, u)|^s f(u)\, du\, \Pi(d\theta) < \infty \qquad \text{for all } s \in (-1, r_0].$$

*Moreover, define*

$$(3.2) \qquad \rho = \exp\left(\int_\Theta \int_\mathbb{R} \log(w(\theta, u)) f(u)\, du\, \Pi(d\theta)\right).$$

*Then for any $\varepsilon > 0$, there exists a bounded function $\nu : \Theta \to \mathbb{R}$, such that*

$$(3.3) \qquad \sup_{\theta \in \Theta} |E(\nu(\theta_1^*) - \nu(\theta) + \log w(\theta, e_1) | \theta_0^* = \theta) - \log \rho| < \varepsilon.$$

Under the stationary distribution $\Pi$, $z(\theta_0^*, e_1)/w(\theta_0^*, e_1) \stackrel{\mathrm{D}}{=} \theta_{0,1}^*$, and therefore $\log \rho$ is also equal to $\int_\Theta \int_\mathbb{R} \log(|z(\theta, u)|/|\theta_1|) f(u)\, du\, \Pi(d\theta)$.

THEOREM 3.2. *Assume Assumptions* A.1–A.6 *and let $\rho$ be as in* (3.2). *If $\rho < 1$, then for any $\rho_1 \in (\rho, 1)$ there exist $K < \infty$, $s > 0$ and $V : \mathbb{R}^p \to [1, \infty)$ such that*

$$E(V(X_1) | X_0 = x) \leq K + \rho_1^s V(x) \qquad \text{for all } x \in \mathbb{R}^p,$$

*and the test function has the form $V(x) = 1 + \lambda(x)\|x\|^s$, where $\lambda$ is bounded and bounded away from 0. As a consequence, $\{X_t\}$ is $V$-uniformly ergodic.*

Although $\rho < 1$ also implies the homogeneous process $\{X_t^*\}$ does not explode, it actually is transient in this case, diminishing to 0.

The value $\log \rho$ turns out to be the Lyapounov exponent for the original chain. This means the condition in Theorem 3.2 is sharp: $\{X_t\}$ is transient if $\rho > 1$ by (3.4) below and Cline and Pu [(2001), Theorems 2.1 and 2.2].

THEOREM 3.3. *Assume Assumptions* A.1–A.6 *and let $\rho$ be as in* (3.2). *Then $\log \rho$ is the Lyapounov exponent for* $\{X_t\}$. *Indeed,*

$$(3.4) \qquad \lim_{n \to \infty} \limsup_{\|x\| \to \infty} \left| \frac{1}{n} E\left( \log\left( \frac{\|X_n\|}{\|X_0\|} \right) \Big| X_0 = x \right) - \log \rho \right| = 0.$$

The key to the piggyback method is the use of the *near-equilibrium equation* (3.3) in the proof of Theorem 3.2. We present here an outline of the method as it is applied to the TAR-ARCH model. This can be used as a guide to the series of lemmas that constitute the proof.

STEP 1. The preliminary step (Lemmas 5.1 and 5.2) is to show that the Markov chains $\{\Theta_t^*\}$ and $\{X_t\}$ are aperiodic, $\phi$-irreducible and sufficiently smooth ($T$-chains), at least when restricted appropriately.

STEP 2. The collapsed process $\{\theta_t^*\}$ stays well within $\Theta_\# = \Theta \setminus H_0$, that is, away from the axial planes, with high probability (Lemma 5.3).

STEP 3. The collapsed process is uniformly ergodic and $\rho$ is well defined (proof of Theorem 3.1).

STEP 4. Let $H_\#$ be the smallest set closed under the map $F\theta = (\theta_2, \ldots, \theta_p, \theta_1)$ and containing the thresholds and axial planes. The collapsed process stays well within $\Theta_{\#\#} = \Theta \setminus H_\#$, with high probability (Lemma 5.4).

STEP 5. Continuity of a function $q(\theta)$ on $\Theta_{\#\#}$ ensures continuity of conditional expectations $E(q(\theta_t^*)|\theta_0^* = \theta)$ on $\Theta_{\#\#}$ (Lemma 5.5). This is required in order to construct a test function that can be piggybacked.

STEP 6. The near-equilibrium equation (3.3) holds with a function $\nu(\theta)$ that is continuous on $\Theta_{\#\#}$ (Lemma 5.6).

STEP 7. Let $\widetilde{X}_t^* = (z(X_{t-1}, e_t), X_{t-1,1}, \ldots, X_{t-1,p-1})$. Both $X_t/\|X_t\|$ and $\widetilde{X}_t^*/\|\widetilde{X}_t^*\|$ are well within $\Theta_{\#\#}$ and close to each other with high probability, if $t$ and $\|X_0\|$ are large enough (Lemmas 5.7 and 5.8). Note that $\widetilde{X}_t^*/\|\widetilde{X}_t^*\|$ is what the first step of the collapsed process would be, if started at $X_{t-1}/\|X_{t-1}\|$.

STEP 8. If, for some function $V_1(x)$ and some $n \geq 1$, $\log(V_1(X_{n+1})/V_1(X_n))$ has negative expectation when conditioned on $X_0 = x$, uniformly for large $\|x\|$, then a test function can be constructed to verify $V$-uniform ergodicity of $\{X_t\}$ (Lemmas 5.9 and 5.10).

STEP 9. The continuity of $\nu(\theta)$ on $\Theta_{\#\#}$, the near-equilibrium equation and the results of Step 7 make it possible to satisfy the condition in Step 8, thus creating a test function for $X_t$ by piggybacking on the simpler properties of the collapsed process (proof of Theorem 3.2).

STEP 10. The Lyapounov exponent for $\{X_t\}$ is $\log \rho$ (proof of Theorem 3.3).

STEP 11. If an appropriate condition on the collapsed process holds, then the above can be applied to a test function roughly proportional to $\|x\|^r$, thus ensuring the stationary distribution of $\{X_t\}$ exists and has $r$th moment (proof of Theorem 3.5 below).

The method just outlined is necessary because of the discontinuities in $a^*$ and $b^*$. One may think that a test function that gives the required drift condition when applied to the homogeneous process $\{X_t^*\}$ could just as well be applied to the original process $\{X_t\}$. Unfortunately this is not the case because the expected difference between $V(X_1)/V(X_0)$ and $V(X_1^*)/V(X_0)$ need not diminish sufficiently fast as $X_0$ increases in magnitude. Nevertheless, as indicated in (3.5) and (3.6) below, solving the simpler problem does suffice to verify stability. This is in fact the primary motivation behind our method.

The piggyback method of construction of a test function depends on obtaining a solution to the near-equilibrium equation (3.3). Indeed, solving (3.3) is tantamount to solving the stability of $\{X_t\}$. It has the advantage that it does not actually require finding the stationary distribution $\Pi$, as the next result indicates, but it does require constructing the function $\nu$.

COROLLARY 3.4. *Assume Assumptions* A.1–A.6. *If there exists a bounded function* $\nu : \Theta \to \mathbb{R}$ *such that*

$$(3.5) \qquad \sup_{\theta \in \Theta} E(\nu(\theta_1^*) - \nu(\theta) + \log w(\theta, e_1) | \theta_0^* = \theta) < 0,$$

*then* $\{X_t\}$ *is $V$-uniformly ergodic.*

Theorem 3.2 does not guarantee the existence of any particular moments of the stationary distribution for $\{X_t\}$. Two equivalent conditions which do are provided next.

THEOREM 3.5. *Assume Assumptions* A.1–A.6 *and suppose* $E(|e_1|^r) < \infty$, $r > 0$. *The following are equivalent conditions for* $\{X_t\}$ *to be* $V$-*uniformly ergodic with* $K\|x\|^r \leq V(x) \leq L + M\|x\|^r$ *for some positive* $K$ *and finite* $L, M$, *and for* $\{X_t\}$ *to have a stationary distribution with finite* $r$*th moment:*

(i) *There exists* $\lambda: \Theta \to (0, \infty)$, *bounded and bounded away from* 0, *such that*

$$\sup_{\theta \in \Theta} E\left(\frac{\lambda(\theta_1^*)}{\lambda(\theta)}(w(\theta, e_1))^r \Big| \theta_0^* = \theta\right) < 1. \tag{3.6}$$

(ii) $\limsup_{n \to \infty} \sup_{\theta \in \Theta} (E(\prod_{t=1}^{n}(w(\theta_{t-1}^*, e_t))^r | \theta_0^* = \theta))^{1/n} < 1.$

We observe that the above conditions are not equivalent to the condition that $E((w(\theta_0^*, e_1))^r) < 1$ under the stationary distribution for $\{\theta_t^*\}$. Though more appealing, and sufficient for verifying ergodicity, the latter is neither necessary nor sufficient for the $r$th moment to exist. (See the examples below.)

For nonthreshold ARCH and GARCH models, the stationary distribution has been shown to have regularly varying tails [e.g., Borkovec (2000), Borkovec and Klüppelberg (2001) and Basrak, Davis and Mikosch (2002)]. (See Example 4.2.) The index of regular variation is the supremum of $r$ satisfying the conditions in Theorem 3.5. Whether this is also true for more general models is as yet an open question.

In general, it seems that only numerical methods will verify the conditions for stability and moments. We recommend simulating $\{\theta_t^*\}$ and estimating the Lyapounov exponent $\log \rho$ in (3.2). There is a distinct advantage, furthermore, to having expressed stability of $\{X_t\}$ in terms of $\{\theta_t^*\}$ because the latter is uniformly ergodic. This means both that convergence of a simulation will tend to be faster and that estimators are well behaved with relatively small variances. The only two alternatives are to conduct the very high-dimensional optimization required to find a test function $V(x)$ or to simulate the time series itself, and pay attention only to its (highly volatile) behavior when it grows large.

**4. Examples.** In this section we provide four examples, giving more specifics about their stability conditions. First, the case $p = 1$ is examined explicitly and conditions given in simple terms. Second, the known criterion for ARCH models is related to the piggyback method, showing how best to estimate the Lyapounov exponent. Third, a special TARCH(2) model can be related to an ARCH(2) model and, finally, second moment conditions for general TARCH($p$) models are considered.

EXAMPLE 4.1 [The TAR-ARCH(1) model]. Assume $p = 1$ and Assumptions A.1 and A.2. Let

$$X_t = a^*(X_{t-1}) + b^*(X_{t-1})e_t + a_0(X_{t-1}) + b_0(X_{t-1})e_t$$

with $a^*(x) = (a_1 \mathbb{1}_{x<0} + a_2 \mathbb{1}_{x>0})x$ and $b^*(x) = (b_1 \mathbb{1}_{x<0} + b_2 \mathbb{1}_{x>0})|x|$, $b_1, b_2 > 0$. Also, $a_0(\cdot)$ and $b_0(\cdot)$ are locally bounded; $a_0(x) = o(\|x\|)$ and $b_0(x) = o(\|x\|)$ as $\|x\| \to \infty$; and $b^*(x) + b_0(x) \neq 0$ for all $x \in \mathbb{R}$. Then Assumptions A.1–A.6 are all satisfied.

Though expressed differently, this includes the usual formulation for a first-order AR-ARCH model with $a_1 = a_2$ and $b_1 = b_2$ [e.g., in Borkovec and Klüppelberg (2001)].

Let $p_1 = P(\theta_1^* = 1|\theta_0^* = -1) = P(a_1 - b_1 e_1 < 0)$ and $p_2 = P(\theta_1^* = -1|\theta_0^* = 1) = P(a_2 + b_2 e_1 < 0)$. Then $\Pi$ is found to be the stationary distribution of a two-state Markov chain, given by $\Pi(-1) = p_2/(p_1 + p_2)$ and $\Pi(1) = p_1/(p_1 + p_2)$. By Theorem 3.2, $\{X_t\}$ is geometrically ergodic if

$$\log \rho = \frac{p_2 E(\log|a_1 - b_1 e_1|) + p_1 E(\log|a_2 + b_2 e_1|)}{p_1 + p_2} < 0.$$

The proof of stability is based on an equilibrium equation (see Lemma 5.6),

$$E(\nu(\theta_1^*) - \nu(\theta) + \log w(\theta, e_1)|\theta_0^* = \theta) = \log \rho, \qquad \theta = \pm 1,$$

which is easily solved here:

$$\nu(\pm 1) = \pm \frac{E(\log|a_2 + b_2 e_1|) - E(\log|a_1 - b_1 e_1|)}{2(p_1 + p_2)}.$$

The function $\lambda$ in (3.6) has a similar, if more cumbersome, solution. To verify that condition, we must find $\gamma$ [$= \lambda(1)/\lambda(-1)$] such that

$$\gamma E(|a_1 - b_1 e_1|^r \mathbb{1}_{a_1 - b_1 e_1 \leq 0}) + E(|a_1 - b_1 e_1|^r \mathbb{1}_{a_1 - b_1 e_1 > 0}) < 1$$

and

$$\gamma^{-1} E(|a_2 + b_2 e_1|^r \mathbb{1}_{a_2 + b_2 e_1 \leq 0}) + E(|a_2 + b_2 e_1|^r \mathbb{1}_{a_2 + b_2 e_1 > 0}) < 1.$$

Letting $E_{i,j} = E(|a_i + (-1)^i b_i e_1|^r \mathbb{1}_{(-1)^j(a_i + (-1)^i b_i e_1) > 0})$, the existence of such a $\gamma$ equates to

(4.1) $\quad \max(E_{1,2}, E_{2,2}) < 1 \quad \text{and} \quad E_{1,1} E_{2,1} < (1 - E_{1,2})(1 - E_{2,2}).$

On the other hand, the condition $E((w(\theta_0^*, e_1))^r) < 1$ under stationarity is

$$\frac{p_2(E_{1,1} + E_{1,2}) + p_1(E_{2,1} + E_{2,2})}{p_1 + p_2} < 1,$$

which certainly neither implies nor is implied by (4.1). Two cases where they do agree are for the ARCH(1) models ($a_1 = a_2 = 0$, $b_1 = b_2$) and for the TARCH(1) models ($a_1 = a_2 = 0$) with errors symmetric about 0.

EXAMPLE 4.2 [The ARCH($p$) model]. Here, the model is

$$\xi_t = (b_0^2 + b_1^2 \xi_{t-1}^2 + \cdots + b_p^2 \xi_{t-p}^2)^{1/2} e_t.$$

Assume Assumptions A.1 and A.2 and each $b_i > 0$. A standard way to handle this is to embed it in a random coefficients model [Bougerol and Picard (1992b) and Basrak, Davis and Mikosch (2002)]. To do this, let $Y_t = (\xi_t^2, \ldots, \xi_{t-p+1}^2)$ so that $Y_t = C_t + B_t Y_{t-1}$ for an i.i.d. sequence of random matrices and vectors $\{(B_t, C_t)\}$. Indeed,

$$B_t = \begin{pmatrix} b_1^2 e_t^2 & b_2^2 e_t^2 & \cdots & b_p^2 e_t^2 \\ 1 & 0 & \cdots & 0 \\ \vdots & \ddots & \ddots & \vdots \\ 0 & \cdots & 1 & 0 \end{pmatrix} \quad \text{and} \quad C_t = \begin{pmatrix} b_0^2 e_t^2 \\ 0 \\ \vdots \\ 0 \end{pmatrix}.$$

Define $M_t = B_t \cdots B_1$, $M_0 = I$ and $\Lambda_t = \log(\|M_t\|/\|M_{t-1}\|)$, where any matrix norm may be chosen. There exists, irrespective of the norm, $\gamma \stackrel{\text{def}}{=} \lim_{t \to \infty} \frac{1}{t} E(\log(\|M_t\|))$ and hence

$$(4.2) \qquad \lim_{t \to \infty} \frac{1}{t} \sum_{i=1}^{t} \Lambda_i = \lim_{t \to \infty} \frac{1}{t} \log(\|M_t\|) = \gamma \qquad \text{a.s.}$$

[Furstenberg and Kesten (1960) and Kingman (1973)]. Under the given assumptions, the necessary and sufficient condition for $\{Y_t\}$ to have a stationary solution is $\gamma < 0$ [Brandt (1986), Bougerol and Picard (1992a), Theorem 2.5, and Goldie and Maller (2000)].

Moreover, suppose $\{\theta_t^*\}$ is the collapsed process for $\{X_t^*\}$ and define $T_t^* = ((\theta_{t,1}^*)^2, \ldots, (\theta_{t,p}^*)^2)'$. Let $\underline{1} = (1, \ldots, 1)'$. Then it is easy to show that

$$T_t^* = \frac{M_t T_0^*}{\underline{1}' M_t T_0^*} \quad \text{for all } t \quad \text{and} \quad (w(\theta_{t-1}^*, e_t))^2 = \frac{\underline{1}' M_t T_0^*}{\underline{1}' M_{t-1} T_0^*}.$$

Therefore, with $T_0^* = \frac{1}{p}\underline{1}$ and $\rho$ as defined in Theorem 3.1,

$$2 \log \rho = \lim_{t \to \infty} \frac{1}{t} \sum_{i=1}^{t} 2 \log(w(\theta_{i-1}^*, e_i))$$

$$= \lim_{t \to \infty} \frac{1}{t} \log(\underline{1}' M_t T_0^*) = \lim_{t \to \infty} \frac{1}{t} \log(\|M_t\|) = \gamma.$$

This gives an alternative method for estimating the Lyapounov exponent. While earlier authors have recommended simulating the matrices $B_t$, computing a norm of their product and applying (4.2) in order to estimate $\gamma$, we recommend simulation of $\theta_t^*$ (or $T_t^*$) instead.

The stationary distribution for related GARCH models is known to have regularly varying tails [Basrak, Davis and Mikosch (2002)]. Their argument applies here for an ARCH($p$) model as well, and the index of regular variation is the positive value of $\kappa$ satisfying

$$(4.3) \qquad \lim_{t \to \infty} (E(\|M_t\|^{\kappa/2}))^{1/t} = 1,$$

which is based on a result of Kesten (1973). See also Goldie (1991). The unbounded support of the error density $f$ (Assumption A.1) ensures that a solution to (4.3) does indeed exist. From our discussion above, we may easily see that such $\kappa$ as satisfies (4.3) must also satisfy

$$\lim_{t \to \infty} \left( E\left( \prod_{i=1}^{t} (w(\theta_{i-1}^*, e_i))^\kappa \right) \right)^{1/t} = 1.$$

Note that the condition in Theorem 3.5(ii) essentially is that $r < \kappa$.

EXAMPLE 4.3 [The TARCH(2) model with delay specific conditional heteroscedasticity].  Let $b_1 > 0$, $b_2 > 0$ and consider the order 2 model:

$$\xi_t = b(X_{t-1})e_t = \begin{cases} (b_{10}^2 + b_1^2(\xi_{t-1}^2 + \xi_{t-2}^2))^{1/2} e_t, & \text{if } \xi_{t-1} \leq 0, \\ (b_{20}^2 + b_2^2(\xi_{t-1}^2 + \xi_{t-2}^2))^{1/2} e_t, & \text{if } \xi_{t-1} > 0. \end{cases}$$

Assume Assumptions A.1 and A.2 and also that the errors have density symmetric about 0. This is a restricted form of the order 2 model with threshold delay 1. Other than the coefficient due to the delay criterion, the conditional heteroscedasticity is proportional to $\|X_t\|$. We will show this model has the same stability criterion as a random coefficients model suggested by the piggyback method, even though the TARCH model itself cannot be embedded in a random coefficients model.

Let $\theta = (\theta_1, \theta_2)$ and $b^*(\theta) = b_1 \mathbb{1}_{\theta_1 \leq 0} + b_2 \mathbb{1}_{\theta_1 > 0}$. Given initial state $\theta_0^* = (\theta_{0,1}^*, \theta_{0,2}^*)$, we have

$$\theta_1^* = (\theta_{1,1}^*, \theta_{1,2}^*) = \frac{(b^*(\theta_0^*)e_1, \theta_{0,1}^*)}{(b^*(\theta_0^*)^2|e_1|^2 + |\theta_{0,1}^*|^2)^{1/2}}.$$

Since the error distribution is symmetric, $|e_1|$, $\text{sgn}(e_1)$ and $\theta_0^*$ are independent. Since $b^*(\theta_0^*)$ depends only on $\text{sgn}(\theta_{0,1}^*) = \text{sgn}(\theta_{1,2}^*)$, it is easy to see then that $\text{sgn}(\theta_{1,1}^*)$ and $(|\theta_{1,1}^*|, \text{sgn}(\theta_{1,2}^*))$ are independent. Thus, under the stationary distribution, $\text{sgn}(\theta_{0,1}^*)$ and $|\theta_{0,1}^*|$ are independent and hence $|\theta_{0,1}^*|$ and $b^*(\theta_0^*)|e_1|$ are independent.

Now consider an ARCH(2) process with uniform heteroscedasticity and error density given by $\bar{f}(u) = \frac{1}{2b_1} f(u/b_1) + \frac{1}{2b_2} f(u/b_2)$. That is, let $\{\bar{e}_t\}$ be an i.i.d. sequence from $\bar{f}$ and $\bar{\xi}_t = (1 + \bar{\xi}_{t-1}^2 + \bar{\xi}_{t-2}^2)^{1/2} \bar{e}_t$. This is, in fact, a special case of the model of this example, but with both coefficients equal to 1 and with a different error distribution now depending on $b_1$ and $b_2$. Applying the comments above accordingly, it is apparent that the stationary distribution of $|\theta_{t,1}^*|$ (but not of $\theta_t^*$) is nevertheless the same for both.

By the comment following Theorem 3.1, the condition for geometric ergodicity of $\{\xi_t\}$ may be expressed as $E(\log(|b^*(\theta_0^*)e_1|/|\theta_{0,1}^*|)) < 0$ under the stationary distribution, or equivalently as $E(\log(|\bar{e}_1|)) - E(\log(|\theta_{0,1}^*|)) < 0$. This depends only on the stationary distribution of $|\theta_{t,1}^*|$. Therefore, $\{\xi_t\}$ and $\{\bar{\xi}_t\}$ have identical stability criteria. The previous example discusses ARCH models in more detail and relates them to a random coefficients model.

The second moment condition for the random coefficients model associated here is $(b_1^2 + b_2^2)E(e_1^2) = E(2\bar{e}_1^2) < 1$, and this coincides with (4.4) found in Example 4.4 for a more general TARCH model.

EXAMPLE 4.4 [The TARCH($p$) model with delay 1]. Assume Assumptions A.1 and A.2. Consider $\xi_t = b(X_{t-1})e_t$, where

$$b(x) = \begin{cases} (b_{10}^2 + b_{11}^2 x_1^2 + \cdots + b_{1p}^2 x_p^2)^{1/2}, & \text{if } x_1 \leq 0, \\ (b_{20}^2 + b_{21}^2 x_1^2 + \cdots + b_{2p}^2 x_p^2)^{1/2}, & \text{if } x_1 > 0, \end{cases}$$

and $b_{ji} > 0$, $i = 1, \ldots, p$, $j = 1, 2$. To get an explicit $r$th moment condition for this model, when $r \leq 2$, define $p_j = P((-1)^j e_1 > 0)$ and $E_j = E(|e_1|^r \mathbb{1}_{(-1)^j e_1 > 0})$, $j = 1, 2$. Then

$$(4.4) \qquad b_{11}^r E_1 + b_{21}^r E_2 + \sum_{i=2}^{p}(b_{1i}^r p_1 + b_{2i}^r p_2)E(|e_1|^r) < 1$$

implies (3.6) is satisfied if $r \leq 2$. This is proven at the end of Section 5. With symmetric errors, this reduces to

$$\tfrac{1}{2}((b_{11}^r + \cdots + b_{1p}^r) + (b_{21}^r + \cdots + b_{2p}^r))E(|e_1|^r) < 1.$$

(Compare this to the condition given in Corollary 2.2.) In fact, when $r = 2$ the test function is optimal, with the expectation in (3.6) having the same value for all $\theta$, suggesting this is the best possible condition for ergodicity with finite second moment. On the other hand, the condition $E((w(\theta_0^*, e_1)^2) < 1$ under stationarity is

$$\int_\Theta ((b_{11}^2 \theta_1^2 + \cdots + b_{1p}^2 \theta_p^2)\mathbb{1}_{\theta_1 \leq 0}$$
$$+ (b_{21}^2 \theta_1^2 + \cdots + b_{2p}^2 \theta_p^2)\mathbb{1}_{\theta_1 > 0} + 1 - \theta_p^2)\Pi(d\theta)E(|e_1|^2) < 1.$$

Simulations indicate that this is not the same as (4.4) with $r = 2$.

**5. Proofs.** The proof for the case $p = 1$ follows the same piggyback principle as does the proof for the case $p > 1$ but it is much simpler, since $\Theta = \{-1, 1\}$ is finite, and we therefore omit it.

We start by asserting some regularity on $\{X_t\}$ and proving the simple conditions in Theorem 2.1 and Corollary 2.2.

LEMMA 5.1. *Assume Assumption* A.1. $\{X_t\}$ *is an aperiodic, $\mu_p$-irreducible $T$-chain on $\mathbb{R}^p$.*

PROOF. The lemma follows from Assumption A.1 by Theorem 2.2(ii) of Cline and Pu (1998). □

PROOF OF THEOREM 2.1. We use essentially the same argument as Lu and Jiang (2001). Let $\beta \in (0,1)$ satisfy $\sum_{i=1}^{p} \beta^{-i} c_i = 1$ and define

$$d_i = \sum_{j=i}^{p} \beta^{i-j-1} c_j, \qquad i = 1, \ldots, p, \ d_{p+1} = 0,$$

so that $\beta d_i = c_i + d_{i+1}$ and $d_1 = 1$. Let $V(x) = 1 + \sum_{i=1}^{p} d_i |x_i|^r$. By (2.2),

$$E(V(X_1)|X_0 = x) = 1 + E(|\xi_1|^r | X_0 = x) + \sum_{i=1}^{p-1} d_{i+1} |x_i|^r$$

$$\leq 1 + K + \sum_{i=1}^{p} c_i |x_i|^r + \sum_{i=1}^{p-1} d_{i+1} |x_i|^r = 1 - \beta + K + \beta V(x).$$

By a standard result [cf. Meyn and Tweedie (1993), Theorem 16.0.1] and Lemma 5.1, $\{X_t\}$ is $V$-uniformly ergodic and the stationary distribution has finite $r$th moment. □

PROOF OF COROLLARY 2.2. Note first that, if $r \leq 2$, then

(5.1) $\quad b^r(x) \leq (b_0^2 + b_1^2 x_1^2 + \cdots + b_p^2 x_p^2)^{r/2} \leq b_0^r + b_1^r |x_1|^r + \cdots + b_p^r |x_p|^r.$

Also, if $r \leq 1$, then

(5.2) $\quad\quad\quad\quad\quad |a(x)|^r \leq a_0^r + a_1^r |x_1|^r + \cdots + a_p^r |x_p|^r.$

However, if $1 < r \leq 2$, then by Jensen's inequality,

(5.3) $\quad (a_1 |x_1| + \cdots + a_p |x_p|)^r \leq (a_1 |x_1|^r + \cdots + a_p |x_p|^r)(a_1 + \cdots + a_p)^{r-1}.$

Also, for any small $\delta > 0$ there exists $M < \infty$ such that $|a(x)|^r \leq M + (1 + \delta)(\sum_{i=1}^{p} a_i |x_i|)^r$.

(i) Since $r \leq 1$, it is immediate that $E(|\xi_1|^r \mid X_0 = x) \leq |a(x)|^r + b^r(x) \times E(|e_1|^r)$, and thus (2.2) follows from (5.1) and (5.2).

(ii) For $1 < r \leq 2$ the conditions on the errors imply

$$E(|\xi_1|^r | X_0 = x) = \tfrac{1}{2} E(|a(x) + b(x)e_1|^r) + \tfrac{1}{2} E(|a(x) - b(x)e_1|^r)$$
$$\leq |a(x)|^r + b^r(x) E(|e_1|^r),$$

using the fact that $(1-u)^r + (1+u)^r - 2u^r \leq 2$ for $0 \leq u \leq 1$. We now apply (5.1) and (5.3) with $c_i = (1+\delta)a_i(\sum_{k=1}^p a_k)^{r-1} + b_i^r E(|e_1|^r)$ and $\delta$ sufficiently small.  $\square$

Next, we have two lemmas about $\{\theta_t^*\}$ followed by the proof of Theorem 3.1.

LEMMA 5.2.   *Assume Assumptions* A.1, A.4 *and* A.5. *Then the following hold:*

(i) $\{X_t^*\}$ *is an aperiodic, $\mu_p$-irreducible $T$-chain when restricted to* $\mathbb{R}_\#^p = \mathbb{R}^p \setminus H_0$.

(ii) $\{\theta_t^*\}$ *is an aperiodic, $\mu_\Theta$-irreducible $T$-chain when restricted to* $\Theta_\# = \Theta \setminus H_0$.

PROOF.   (i) By Assumption A.5, $b^*$ is locally bounded away from 0 on $\mathbb{R}_\#^p$, so the $n$-step transition density $g_n^\#(\cdot; x)$ is well defined on $\mathbb{R}_\#^p$ for any $n \geq p$ and any initial state $x$ in $\mathbb{R}_\#^p$. Furthermore, it is locally bounded away from 0 and full. Hence $\mathbb{R}_\#^p$ is absorbing and $\{X_t^*\}$ is aperiodic and $\mu_p$-irreducible on $\mathbb{R}_\#^p$.

For any $z \in \mathbb{R}_\#^p$, let $G_z$ be an open set containing $z$ such that its closure is compact and contained in $\mathbb{R}_\#^p$. Then

$$T_z(x, A) \stackrel{\text{def}}{=} \mathbb{1}_{G_z}(x) \int_A \inf_{u \in G_z} g_p^\#(y; u) \, dy$$

defines a kernel on $\mathbb{R}_\#^p$, nontrivial at $z$, such that $T_z(\cdot, A)$ is lower semicontinuous for each $A$ and $T_z(x, A) \leq P^p(x, A)$. By Meyn and Tweedie [(1993), Proposition 6.2.4], $\{X_t^*\}$ is a $T$-chain when restricted to $\mathbb{R}_\#^p$.

(ii) For any set $A \subset \Theta$ we define the cone $C_A = \{x \in \mathbb{R}_0^p : x/\|x\| \in A\}$. Since

$$P(\theta_t^* \in A | \theta_0^* = \theta) = P(X_t^* \in C_A | X_0^* = c\theta)$$

for any $A$, $t$, $c$ and $\theta$, the result follows easily from (i).  $\square$

LEMMA 5.3.   *Assume Assumptions* A.1–A.5. *Given* $\varepsilon_1 \in (0,1)$, *there exists compact* $C_1 \subset \Theta_\#$ *such that the following hold:*

(i) $P(\theta_1^* \in C_1 | \theta_0^* = \theta) > 1 - \varepsilon_1$ *for all* $\theta \in C_1$;
(ii) $P(\theta_p^* \in C_1 | \theta_0^* = \theta) > 1 - p\varepsilon_1$ *for all* $\theta \in \Theta$.

PROOF. Let $L_0 = \sup_u (1+|u|) f(u)$, which is finite by Assumption A.2. Suppose $\alpha$ and $\beta$ are values with $\beta \geq 0$ and $\max(|\alpha|, \beta) > 0$, and suppose $\varepsilon \in (0, \max(|\alpha|, \beta)/2)$. This implies $\beta + \max(|\alpha| - \varepsilon, 0) > \max(|\alpha|, \beta)/2$. Then, if $\beta > 0$,

$$P(|\alpha + \beta e_t| \leq \varepsilon) \leq L_0 \int_{(-\alpha-\varepsilon)/\beta}^{(-\alpha+\varepsilon)/\beta} \frac{1}{1+|u|} \, du$$

(5.4)

$$\leq \frac{2L_0 \varepsilon}{\beta + \max(|\alpha| - \varepsilon, 0)} \leq \frac{4L_0 \varepsilon}{\max(|\alpha|, \beta)}.$$

The above holds trivially in case $\beta = 0$ since, in that case, $|\alpha + \beta e_t| = |\alpha| > \varepsilon$.

Now let $L_1 = \inf_{\theta \in \Theta} \max(|a^*(\theta)|, b^*(\theta))$ and $L_2 = \sup_{\theta \in \Theta} \max(|a^*(\theta)|, b^*(\theta))$. Note $L_1 > 0$ by Assumption A.5 and $L_2 < \infty$ by Assumption A.3. By (5.4), we determine that if $\varepsilon \in (0, L_1/2)$, then

(5.5) $$P(|z(\theta, e_t)| \leq \varepsilon) \leq 4L_0 \varepsilon / L_1.$$

With no loss we suppose $\varepsilon_1 \in (0, L_1/2)$. Let $\varepsilon_2 = \frac{\varepsilon_1}{1 + 4L_0/L_1}$. Choose $M_1 \geq 1$ so that $P(|e_t| > M_1) < \varepsilon_2$ and thus, applying (5.5) with $\varepsilon = \varepsilon_2$,

(5.6) $$P(|z(\theta, e_t)| > \varepsilon_2; |e_t| \leq M_1) > 1 - \varepsilon_1.$$

Fix $\gamma = 1 + L_2(1 + M_1)$ and define $X_{1,x,u}^* = (z(x,u), x_1, \ldots, x_{p-1})$. Then

(5.7) $\quad |u| \leq M_1 \implies |z(\theta, u)| \leq \gamma - 1 \quad \text{and} \quad \|X_{1,\theta,u}^*\| \leq |z(\theta, u)| + 1 \leq \gamma.$

By (5.6) and (5.7),

(5.8)
$$P(|z(\theta, e_1)| \geq \varepsilon_2; \|X_{1,\theta,e_1}^*\| \leq \gamma)$$
$$\geq P(|z(\theta, e_1)| \geq \varepsilon_2; |e_1| \leq M_1) > 1 - \varepsilon_1.$$

Now define $C_1 = \{\theta = (\theta_1, \ldots, \theta_p) : |\theta_i| \geq \varepsilon_2 \gamma^{-i}, i = 1, \ldots, p\}$, which is compact. If $\theta \in C_1$, $|z(\theta, e_1)| \geq \varepsilon_2$ and $\|X_{1,\theta,e_1}^*\| \leq \gamma$, then $|z(\theta, e_1)|/\|X_{1,\theta,e_1}^*\| \geq \varepsilon_2 \gamma^{-1}$ and $|\theta_i|/\|X_{1,\theta,e_1}^*\| \geq \varepsilon_2 \gamma^{-i-1}$ for $i = 1, \ldots, p-1$. Since $\theta_1^* = X_{1,\theta,e_1}^*/\|X_{1,\theta,e_1}^*\|$ when $\theta_0^* = \theta$, therefore, we conclude

$$P(\theta_1^* \in C_1 | \theta_0^* = \theta) \geq P(|z(\theta, e_1)| \geq \varepsilon_2; \|X_{1,\theta,e_1}^*\| \leq \gamma)$$
$$> 1 - \varepsilon_1 \quad \text{for } \theta \in C_1.$$

Thus (i) is proved. By (5.8) and an induction argument, (ii) follows as well. □

PROOF OF THEOREM 3.1. By Lemma 5.2(ii), $\{\theta_t^*\}$ restricted to $\Theta_\#$ is an aperiodic, $\mu_\Theta$-irreducible $T$-chain. For any $\varepsilon_1 \in (0,1)$, let $C_1$ be as in Lemma 5.3. Then $C_1$, being compact, is small for the restricted process. It follows easily that $C_1$ is small for the process on $\Theta$ and therefore

Lemma 5.3(ii) implies $\Theta$ is itself small [Meyn and Tweedie (1993), Proposition 5.5.4(i)]. It also follows from the lemmas that the unrestricted process is $\mu_\Theta$-irreducible and aperiodic. Thus, the chain is uniformly ergodic on $\Theta$ by Meyn and Tweedie [(1993), Theorem 16.2.2] is a $T$-chain and has some stationary distribution $\Pi$.

The inequality (3.1) clearly holds for $s \in [0, r_0]$ by Assumption A.2 and the boundedness of $a^*$ and $b^*$ on $\Theta$. By (5.5), there exists $K < \infty$ such that

$$\int_\Theta \int_\mathbb{R} \mathbb{1}_{|z(\theta,u)| \leq \varepsilon} f(u) \, du \, \Pi(d\theta) < K\varepsilon$$

for any $\varepsilon > 0$. Thus, (3.1) holds for $s \in (-1, 0)$. We note further that $|z(\theta, e_1)| \leq w(\theta, e_1) \leq |z(\theta, e_1)| + 1$. Hence Assumption A.2 and (5.5) imply $\{|\log(w(\theta, e_1))|\}_{\theta \in \Theta}$ is uniformly integrable. Thus $\log(w(\theta_0^*, e_1))$ has finite mean under $\Pi$ and $\rho$ is well defined.

The near-equilibrium equation (3.3) holds by Lemma 5.6 below. $\square$

We now provide several lemmas needed to construct a proof of Theorem 3.2. Recall that by Assumption A.6, $a^*$ and $b^*$ are continuous off a collection of hyperplanes, $H_1, \ldots, H_m$, each of which contains the origin. Define $F$ by $F\theta = (\theta_2, \ldots, \theta_p, \theta_1)$ and $\Theta_{\#\#} = \{\theta \in \Theta_\# : F^k \theta \notin H_j \text{ for all } j, k\}$. Note that $\Theta_{\#\#}$ excludes not only the thresholds but also certain critical points for which the probability of being near the thresholds in the first $p$ steps cannot be controlled.

LEMMA 5.4. *Assume Assumptions* A.1–A.6. *Given any $\varepsilon_3 \in (0, 1)$, there exists compact $D_1 \subset \Theta_{\#\#}$ such that*

$$P(\theta_1^* \in D_1 | \theta_0^* = \theta) > 1 - \varepsilon_3 \qquad \text{for all } \theta \in D_1.$$

PROOF. There exist $h_1, \ldots, h_m \in \mathbb{R}_0^p$ such that $\theta \in \Theta_{\#\#}$ if and only if $h_j' F^k \theta \neq 0$ for all $j, k$. Since $F^p = I$, it suffices to assume the first coordinate of each $h_j$ is 1.

Set $\varepsilon_1 = \varepsilon_3/2$ and define $\gamma$, $C_1$ and $X_{1,\theta,e_1}^*$ as in the proof of Lemma 5.3. Then, by that proof,

(5.9) $\quad P(\theta_1^* \in C_1, \|X_{1,\theta,e_1}^*\| \leq \gamma | \theta_0^* = \theta) > 1 - \varepsilon_3/2 \qquad$ for all $\theta \in C_1$.

Also let $L_0$ be as in the proof of Lemma 5.3, $L_3 = \inf_{\theta \in C_1} b^*(\theta)$, $L_4 = \min\{|h_{j,k}| : h_{j,k} \neq 0\}$ and $\varepsilon_4 = \frac{L_3 L_4 \varepsilon_3}{4mp\gamma L_0}$. Define

$$D_1 = \{\theta \in C_1 : |h_j' F^k \theta| \geq \varepsilon_4 \gamma^{-k} \text{ for all } j = 1, \ldots, m, k = 0, \ldots, p-1\}.$$

Define $h_{j,p+1} = h_{j,1} = 1$. If $\theta_0^* = \theta$ and $\|X_{1,\theta,e_1}^*\| \leq \gamma$, then

(5.10) $\begin{aligned}|h_j' F^k \theta_1^*| &= |h_{j,p+1-k}(z(\theta, e_1) - \theta_p) + h_j' F^{k-1} \theta| / \|X_{1,\theta,e_1}^*\| \\ &\geq |h_{j,p+1-k}(z(\theta, e_1) - \theta_p) + h_j' F^{k-1} \theta| / \gamma.\end{aligned}$

Thus $\theta \in D_1$ and $h_{j,p+1-k} = 0$ imply

(5.11) $$|h_j' F^k \theta_1^*| \geq |h_j' F^{k-1} \theta|/\gamma \geq \varepsilon_4 \gamma^{-k}.$$

Also, note that, if $\theta \in C_1$, then

(5.12) $$P(|z(\theta, e_1) - u| \leq c) = P(|a^*(\theta) - b^*(\theta)e_1 - u| \leq c) \leq \frac{2L_0}{L_3} c,$$

for all $u \in \mathbb{R}$ and $c > 0$. Let

$$u_{j,k,\theta} = \begin{cases} \theta_p - \dfrac{h_j' F^{k-1} \theta}{h_{j,p+1-k}\gamma}, & \text{if } h_{j,p+1-k} \neq 0, \\ 0, & \text{if } h_{j,p+1-k} = 0. \end{cases}$$

Therefore, using (5.9)–(5.12),

$$P(\theta_1^* \in D_1 | \theta_0^* = \theta)$$
$$\geq P(\theta_1^* \in C_1; \|X_{1,\theta,e_1}^*\| \leq \gamma; |h_j' F^k \theta_1^*| \geq \varepsilon_4 \gamma^{-k} \text{ for all } j,k | \theta_0^* = \theta)$$
$$> 1 - \frac{\varepsilon_3}{2} - P\left(|z(\theta, e_1) - u_{j,k,\theta}| \leq \frac{\varepsilon_\gamma^{1-k}}{L_4} \text{ for some } j,k\right)$$
$$> 1 - \frac{\varepsilon_3}{2} - \frac{2mpL_0 \varepsilon_4 \gamma}{L_3 L_4} = 1 - \varepsilon_3. \qquad \square$$

LEMMA 5.5. *Assume Assumptions* A.1–A.6. *Suppose $q(\theta)$ is bounded on $\Theta$ and continuous on $\Theta_{\#\#}$. Then $E(q(\theta_t^*)|\theta_0^* = \theta)$ is continuous on $\Theta_{\#\#}$ for each $t \geq 1$.*

PROOF. Fix $t \geq 1$, $\varepsilon \in (0,1)$ and compact $C \subset \Theta_{\#\#}$. Define $L_5 = \sup_{\theta \in \Theta} |q(\theta)|$. Pick $\varepsilon_3 \in (0, \frac{\varepsilon}{1+4tL_5})$. By Lemma 5.4 and its proof, which in turn depends on the proof of Lemma 5.3, there exists $M_1$ so large that $P(|e_1| > M_1) < \varepsilon_3/2$ and compact $D_1$ with $C \subset D_1 \subset \Theta_{\#\#}$ such that

(5.13) $$P(\theta_1^* \in D_1 | \theta_0^* = \theta) > 1 - \varepsilon_3 \qquad \text{for all } \theta \in D_1.$$

Again let $X_{1,x,u}^* = (z(x,u), x_1, \ldots, x_{p-1})$ and define $\theta_{1,x,u}^* = X_{1,x,u}^*/\|X_{1,x,u}^*\|$.

Since $q$ is uniformly continuous on $D_1$, choose $\delta_0 > 0$ such that

(5.14) $$\|\theta' - \theta''\| < \delta_0, \qquad \theta' \in D_1,\ \theta'' \in D_1 \quad \Longrightarrow \quad |q(\theta') - q(\theta'')| < \varepsilon_3.$$

Also, $\{\theta_{1,\theta',u}^* : |u| \leq M_1\}$ is uniformly equicontinuous on $D_1$, so choose $\delta_k$, $k \geq 1$, such that if $|u| \leq M_1$,

(5.15) $$\|\theta' - \theta''\| < \delta_k, \qquad \theta' \in D_1, \theta'' \in D_1$$
$$\Longrightarrow \quad \|\theta_{1,\theta',u}^* - \theta_{1,\theta'',u}^*\| < \delta_{k-1}.$$

Now let $\{\theta^*_{k,\theta'}\}$ and $\{\theta^*_{k,\theta''}\}$ be the processes that start with $\theta'$ and $\theta''$, respectively. We thus have, by (5.13) and (5.15),

$$P(\|\theta^*_{k,\theta'} - \theta^*_{k,\theta''}\| < \delta_{t-k}; \theta^*_{k,\theta'} \in D_1; \theta^*_{k,\theta''} \in D_1)$$
$$> P(\|\theta^*_{k-1,\theta'} - \theta^*_{k-1,\theta''}\| < \delta_{t-k+1};$$
$$\theta^*_{k-1,\theta'} \in D_1; \theta^*_{k-1,\theta''} \in D_1; |e_k| \le M_1) - \varepsilon_3$$
$$> P(\|\theta^*_{k-1,\theta'} - \theta^*_{k-1,\theta''}\| < \delta_{t-k+1}; \theta^*_{k-1,\theta'} \in D_1; \theta^*_{k-1,\theta''} \in D_1) - 2\varepsilon_3,$$

for $k = 1, \ldots, t$. Hence, if $\theta' \in D_1$, $\theta'' \in D_1$ and $\|\theta' - \theta''\| < \delta_t$, then

$$P(\|\theta^*_{t,\theta'} - \theta^*_{t,\theta''}\| < \delta_0; \theta^*_{t,\theta'} \in D_1; \theta^*_{t,\theta''} \in D_1) > 1 - 2t\varepsilon_3.$$

Applying (5.14),

(5.16) $$P(|q(\theta^*_{t,\theta'}) - q(\theta^*_{t,\theta''})| < \varepsilon_3) > 1 - 2t\varepsilon_3.$$

From (5.16) we obtain, if $\theta' \in D_1$, $\theta'' \in D_1$ and $\|\theta' - \theta''\| < \delta_t$, then

$$|E(q(\theta^*_t)|\theta^*_0 = \theta') - E(q(\theta^*_t)|\theta^*_0 = \theta'')| \le E(|q(\theta^*_{t,\theta'}) - q(\theta^*_{t,\theta''})|)$$
$$< \varepsilon_3 + 2L_5 P(|q(\theta^*_{t,\theta'}) - q(\theta^*_{t,\theta''})| \ge \varepsilon_3)$$
$$< (1 + 4tL_5)\varepsilon_3 < \varepsilon.$$

The conclusion then follows since $C \subset D_1$ and since $\varepsilon$ and $C$ are arbitrary. □

This lemma identifies the implicit behavior of $\{\theta^*_t\}$ that we will piggyback upon, namely a near-equilibrium equation.

LEMMA 5.6. *Assume Assumptions* A.1–A.6. *If $\rho$ is defined as in* (3.2), *then for any $\varepsilon_5 > 0$, there exists a bounded function $\nu : \Theta \to \mathbb{R}$, such that*

(5.17) $$\sup_{\theta \in \Theta} |E(\nu(\theta^*_1) - \nu(\theta) + \log w(\theta, e_1)|\theta^*_0 = \theta) - \log \rho| < \varepsilon_5.$$

*Moreover, $\nu$ is continuous on $\Theta_{\#\#}$.*

PROOF. Fix $\varepsilon_5 > 0$. By definition,

(5.18) $$\log \rho = \int_\Theta \int_\mathbb{R} \log(w(\theta, u)) f(u) \, du \, \Pi(d\theta).$$

Define $q(\theta) = E(\log w(\theta, e_1))$, which is bounded. By (5.18) and the uniform ergodicity of $\{\theta^*_t\}$, there exists $\delta < 1$ and $K_1 < \infty$ such that

$$|E(q(\theta^*_t)|\theta^*_0 = \theta) - \log \rho| < K_1 \delta^t \qquad \text{for all } t \ge 1 \text{ and all } \theta \in \Theta$$

[cf. Meyn and Tweedie [(1993), Theorem 16.2.1]. Choose $T$ such that $K_1 \delta^T \leq \varepsilon_5$ and let

$$\nu(\theta) = \sum_{t=0}^{T-1} E(q(\theta_t^*)|\theta_0^* = \theta).$$

Then

$$E(\nu(\theta_1^*) - \nu(\theta) + \log w(\theta, e_1)|\theta_0^* = \theta) = E(q(\theta_T^*)|\theta_0^* = \theta)$$

and (5.17) follows immediately.

Furthermore, $q$ is continuous on $\Theta_{\#\#}$ and therefore $\nu$ is also continuous on $\Theta_{\#\#}$, by Lemma 5.5. $\square$

Previously we defined $X_{1,x,u}^* = (z(x,u), x_1, \ldots, x_{p-1})$ and $\theta_{1,x,u}^* = X_{1,x,u}^*/\|X_{1,x,u}^*\|$. Now define $X_{1,x,u} = (a(x) + b(x)u, x_1, \ldots, x_{p-1})$ and $\theta_{1,x,u} = X_{1,x,u}/\|X_{1,x,u}\|$. Additionally, let

$$\tilde{\theta}_t = \theta_{1,X_{t-1},e_t} = X_t/\|X_t\|,$$
$$\widetilde{X}_t^* = X_{1,X_{t-1},e_t}^* \quad \text{and} \quad \tilde{\theta}_t^* = \theta_{1,X_{t-1},e_t}^* = \widetilde{X}_t^*/\|\widetilde{X}_t^*\|.$$

In the following arguments it will be crucial to compare $X_t$ to $\widetilde{X}_t^*$ (and $\tilde{\theta}_t$ to $\tilde{\theta}_t^*$), to consider how they would differ if $X_{t-1}$ is acted upon by $a^*$ and $b^*$ instead of by $a$ and $b$.

Recall $z(x, u) = a^*(x) + b^*(x)u$. Define also $c(x, u) = a_0(x) + b_0(x)u$ so that $a(x) + b(x)u = z(x, u) + c(x, u)$.

The next two lemmas show that $\tilde{\theta}_t$ and $\tilde{\theta}_t^*$ are each eventually within a compact subset of $\Theta_{\#\#}$ (hence away from singular behavior), with high probability.

LEMMA 5.7. *Assume Assumptions* A.1–A.5. *Given any $\varepsilon_1 \in (0, 1)$, there exist $\gamma > 1$, $\eta \in (0, 1)$, $M_2 < \infty$ and compact set $C_2 \subset \Theta_\#$ such that, if $\|x\| > M_2 \eta^{1-t}$ and if either $x/\|x\| \in C_2$ or $t \geq p$, then*

$$(5.19) \quad P\left(\tilde{\theta}_t \in C_2;\ \tilde{\theta}_t^* \in C_2;\ \eta^t \leq \frac{\|X_t\|}{\|x\|} \leq (\gamma + \eta)^t \Big| X_0 = x\right) > (1 - \varepsilon_1)^t.$$

PROOF. As in the proof of Lemma 5.3, we assume $\varepsilon_1 \in (0, L_0/2)$. Let $\varepsilon_2$, $M_1$ and $\gamma$ be the same as in the proof of Lemma 5.3. We set $\eta = \varepsilon_2/2$. By Assumption A.4 and the definitions of $a_0$ and $b_0$ that follow it, there exists $M_2 < \infty$ such that

$$(5.20) \qquad \frac{|c(x, u)|}{\|x\|} < \eta \qquad \text{for all } |u| \leq M_1,\ \|x\| > M_2.$$

For $t \geq 1$ and $x \in \mathbb{R}^p$, $x \neq 0$, let

$$C_{2,x,t} = \left\{\theta : |\theta_i| \geq \eta(\gamma + \eta)^{-i}, i \leq \min(t, p); |\theta_i| \geq \frac{|x_{i-t}|}{\|x\|}(\gamma + \eta)^{-1}, t < i \leq p\right\}.$$

Define $C_2 = C_{2,x,p}$ which in fact does not depend on $x$. Clearly, if $x/\|x\| \in C_2$ or $t \geq p$, then $C_{2,x,t} = C_2$.

Therefore it suffices to show that, for $t \geq 1$ and $\|x\| > M_2/\eta^{t-1}$,

$$(5.21) \quad P\left(\tilde{\theta}_t \in C_{2,x,t}; \tilde{\theta}_t^* \in C_{2,x,t}; \eta^t \leq \frac{\|X_t\|}{\|x\|} \leq (\gamma + \eta)^t \Big| X_0 = x\right) > (1 - \varepsilon_1)^t.$$

Suppose $\|x\| > M_2$, $|u| \leq M_1$ and $|z(x/\|x\|, u)| \geq \varepsilon_2 = 2\eta$. Then, by (5.7) and (5.20),

$$\eta\|x\| \leq |z(x, u) + c(x, u)| \leq \|X_{1,x,u}\|$$
$$\leq |z(x, u) + c(x, u)| + \|x\| \leq (\gamma + \eta)\|x\|.$$

Likewise,

$$2\eta\|x\| \leq \|X_{1,x,u}^*\| \leq \gamma\|x\|.$$

Furthermore, for the components of $\theta_{1,x,u}$ and $\theta_{1,x,u}^*$,

$$\min(|\theta_{1,i,x,u}|, |\theta_{1,i,x,u}^*|) \geq \begin{cases} \eta(\gamma + \eta)^{-1}, & i = 1, \\ \dfrac{|x_{i-1}|}{\|x\|}(\gamma + \eta)^{-1}, & i = 2, \ldots, p. \end{cases}$$

Hence

$$(5.22) \quad \begin{aligned} P&\left(\theta_{1,x,u} \in C_{2,x,1}; \theta_{1,x,u}^* \in C_{2,x,1}; \eta \leq \frac{\|X_1\|}{\|x\|} \leq (\gamma + \eta) \Big| X_0 = x\right) \\ &\geq P\left(|e_1| \leq M_1; \left|z\left(\frac{x}{\|x\|}, e_1\right)\right| > \varepsilon_2\right) > 1 - \varepsilon_1, \end{aligned}$$

where the last inequality follows from (5.6). That is, (5.21) holds for $t = 1$.

Next, assume (5.21) holds for $t = k$, some $k \geq 1$. Then for $\|x\| > M_2/\eta^k$,

$$P\left(\tilde{\theta}_{k+1} \in C_{2,x,k+1}; \tilde{\theta}_{k+1}^* \in C_{2,x,k+1}; \eta^{k+1} \leq \frac{\|X_{k+1}\|}{\|x\|} \leq (\gamma + \eta)^{k+1} \Big| X_0 = x\right)$$

$$\geq P\left(\tilde{\theta}_{k+1} \in C_{2,X_k,1}; \tilde{\theta}_{k+1}^* \in C_{2,X_k,1}; \tilde{\theta}_k \in C_{2,x,k}; \eta \leq \frac{\|X_{k+1}\|}{\|X_k\|} \leq (\gamma + \eta);\right.$$
$$\left. \text{and } \eta^k \leq \frac{\|X_k\|}{\|x\|} \leq (\gamma + \eta)^k \Big| X_0 = x\right)$$

$$\geq (1 - \varepsilon_1)^{k+1},$$

by conditioning on $X_k$. Thus (5.21) holds for all $t \geq 1$ by induction. $\square$

LEMMA 5.8. *Assume Assumptions* A.1–A.6. *Given any $\varepsilon_3 \in (0,1)$, there exist $\gamma > 1$, $\eta \in (0,1)$, $M_2 < \infty$ and compact set $D_2 \subset \Theta_{\#\#}$ such that, if $\|x\| > M_2 \eta^{1-t}$ and if either $x/\|x\| \in D_2$ or $t \geq 2p$, then*

$$P\left(\tilde{\theta}_t \in D_2; \tilde{\theta}_t^* \in D_2; \eta^t \leq \frac{\|X_t\|}{\|x\|} \leq (\gamma+\eta)^t \Big| X_0 = x\right) > (1-\varepsilon_3)^t.$$

PROOF. The argument for this is analogous to that of Lemma 5.7: an induction by conditioning, starting with $t = p$ and (5.19) and iterating to $t \geq 2p$. Using the same notation as in Lemmas 5.4 and 5.7, the induction step is possible by showing that if $\|x\| > M_2 \eta^{-t}$, $\tilde{\theta}_0 = x/\|x\| \in C_2$ and $h'_j F^k \tilde{\theta}_0 \geq \varepsilon_4 (\gamma + \eta)^{-k}$ for $k < t$, $t = 0, \ldots, p-1$, then

$$P(h'_j F^k \tilde{\theta}_1 \geq \varepsilon_4 (\gamma+\eta)^{-k};$$

$$h'_j F^k \tilde{\theta}_1^* \geq \varepsilon_4 (\gamma+\eta)^{-k}, 1 \leq j \leq m, 0 \leq k \leq t; A_1) > 1 - \varepsilon_3/p,$$

where

$$A_1 = \{\tilde{\theta}_1 \in C_2; \tilde{\theta}_1^* \in C_2; \eta\|x\| \leq \|X_1\| \leq (\gamma+\eta)\|x\|\}.$$

In other words, if $\tilde{\theta}_t$ and $\tilde{\theta}_t^*$ are both in $C_2$, then we can reach $D_2$ with both in $p$ steps, with high probability. The argument for this is the same that yielded (5.21) and (5.22). □

LEMMA 5.9. *Assume Assumptions* A.1–A.5. *Then*

$$\sup_{x \in \mathbb{R}^p} E\left(\left(\frac{1+\|X_{1,x,e_1}\|}{1+\|x\|}\right)^s\right) < \infty \qquad \text{for all } s \in (-1, r_0],$$

*and $\|x\| \to \infty$ implies $\|X_{1,x,e_1}\| \to \infty$ in probability.*

PROOF. Let $K_2 = \sup_{x \in \mathbb{R}^p} \frac{\max(|a(x)|, b(x))}{1+\|x\|}$. Clearly, $\frac{1+\|X_{1,x,e_1}\|}{1+\|x\|} \leq 1 + K_2(1+|e_1|)$, which establishes the result for $0 \leq s \leq r_0$. (Assumptions A.4 and A.5 are not needed for this.)

Next, let $L_1 = \inf_{\theta \in \Theta} \max(|a^*(\theta)|, b^*(\theta))$, as in the proof of Lemma 5.3. Choose $M_3 \geq 1$ such that

$$\max\left(\frac{|a_0(x)|}{1+\|x\|}, \frac{|b_0(x)|}{1+\|x\|}\right) \leq \frac{L_1}{4} \qquad \text{whenever } \|x\| > M_3,$$

by the definitions of $a_0(x)$ and $b_0(x)$ following Assumption A.4. Define $L_6 = \inf_{\|x\| \leq M_3} \frac{b(x)}{1+\|x\|}$, which is positive by Assumption A.1. Thus, $\|x\| \leq M_3$ implies

$$\max\left(\frac{|a(x)|}{1+\|x\|}, \frac{b(x)}{1+\|x\|}\right) \geq L_6.$$

On the other hand, $\|x\| > M_3$ implies

$$\max\left(\frac{|a(x)|}{1+\|x\|}, \frac{b(x)}{1+\|x\|}\right) \geq \max\left(\frac{|a^*(x)| - |a_0(x)|}{1+\|x\|}, \frac{b^*(x) - |b_0(x)|}{1+\|x\|}\right)$$

$$\geq \frac{\max(|a^*(x)|, b^*(x))}{2\|x\|} - \frac{L_1}{4} \geq \frac{L_1}{4}.$$

Applying (5.4), then, for all $\varepsilon \in (0, \min(L_1/8, L_6/2))$ and $x \in \mathbb{R}^p$,

$$P\left(\frac{1+\|X_{1,x,e_1}\|}{1+\|x\|} \leq \varepsilon\right) \leq P\left(\left|\frac{a(x)}{1+\|x\|} + \frac{b(x)}{1+\|x\|}e_1\right| \leq \varepsilon\right) \leq \frac{4L_0\varepsilon}{\min(L_1/4, L_6)}.$$

This ensures the expectation in the lemma statement is uniformly bounded for each $s \in (-1, 0)$. Furthermore, it also requires that $\|X_{1,x,e_1}\| \to \infty$ in probability, as $\|x\| \to \infty$. $\square$

LEMMA 5.10. *Assume Assumptions* A.1–A.3. *Suppose there exist a bounded function $\nu(x)$, $\rho_1 < 1$ and $n \geq 0$ such that*

$$\limsup_{\|x\| \to \infty} E\left(\nu(X_{n+1}) - \nu(X_n) + \log\left(\frac{1+\|X_{n+1}\|}{1+\|X_n\|}\right)\bigg| X_0 = x\right) < \log \rho_1 < 0.$$

*Then, for some $s \in (0, r_0]$ and some function $\lambda \colon \mathbb{R}^p \to (0, \infty)$, bounded and bounded away from 0,*

$$(5.23) \qquad \limsup_{\|x\| \to \infty} E\left(\frac{\lambda(X_1)\|X_1\|^s}{\lambda(x)\|x\|^s}\bigg| X_0 = x\right) < \rho_1^s < 1.$$

*Consequently, $\{X_t\}$ is $V$-uniformly ergodic with $V(x) = 1 + \lambda(x)\|x\|^s$ and (2.1) holding.*

PROOF. The proof of (5.23) is contained in the proof of Lemma 4.1(i) of Cline and Pu (1999), with Lemma 5.9 in support. Along with the conclusion of Lemma 5.1 and the fact $E(\lambda(X_1)\|X_1\|^s | X_0 = x)$ is locally bounded, this establishes the standard Foster–Lyapounov drift condition for geometric ergodicity [Meyn and Tweedie (1993), Theorem 15.0.1 and Section 15.2]. $\square$

PROOF OF THEOREM 3.2. We do not assume $\rho < 1$ until the very end of the proof. Choose $\varepsilon_6 > 0$ and let $\varepsilon_5 = \varepsilon_6/4$. Define $\nu$ according to Lemma 5.6. Extend $\nu$ to $\mathbb{R}^p$ by defining $\nu(0) = 0$ and $\nu(x) = \nu(x/\|x\|)$ for $x \neq 0$. Let $V_1(x) = 1 + \|x\|$. Then

$$\limsup_{\|y\| \to \infty} |E(\nu(X^*_{1,y,e_1}) - \nu(y) + \log(V_1(X^*_{1,y,e_1})/V_1(y))) - \log \rho|$$
(5.24)
$$< \varepsilon_5 = \varepsilon_6/4.$$

Let $L_7 = \sup_{x \in \mathbb{R}^p} |\nu(x)|$. By Lemma 5.9, there exists $K_4 < \infty$ such that
$$\sup_{y \in \mathbb{R}^p} E(|\log(V_1(X_{1,y,e_1})/V_1(y))|) \leq K_4.$$

Let $L_8 = 4 \max(2p|\log \rho|, 2(2p+1)L_7, 2p(2L_7 + K_4))$ and choose $\varepsilon_3 < \varepsilon_6/L_8$. Choose $D_2$, $\eta$ and $M_2$ according to Lemma 5.8 and define
$$A_t = \{\tilde{\theta}_t \in D_2; \tilde{\theta}_t^* \in D_2; \|X_t\| \geq \eta^t \|X_0\|\}, \qquad t = 1, \ldots, 2p+1.$$

Then, in particular,

(5.25) $$P(A_t | X_0 = x) > 1 - t\varepsilon_3,$$

if $\|x\| > M_2 \eta^{1-t}$ and either $t \geq 2p$ or $x/\|x\| \in D_2$. Conditioning on $(X_{2p-1}, e_{2p})$ and applying (5.24) and (5.25), we have

(5.26)
$$\limsup_{\|x\| \to \infty} |E((\nu(\widetilde{X}_{2p+1}^*) - \nu(X_{2p}) + \log(V_1(\widetilde{X}_{2p+1}^*)/V_1(X_{2p}))) \mathbb{1}_{A_{2p}} | X_0 = x) - \log \rho|$$
$$< \varepsilon_6/4 + 2p|\log \rho|\varepsilon_3 < \varepsilon_6/2.$$

Since $a_0(x) = o(\|x\|)$, $b_0(x) = o(\|x\|)$ and $D_2 \subset \Theta_\#$, we conclude
$$\liminf_{\substack{\|y\| \to \infty \\ y/\|y\| \in D_2}} \frac{\|X_{1,y,e_1}\|}{\|y\|} = \liminf_{\substack{\|y\| \to \infty \\ y/\|y\| \in D_2}} \frac{\|X_{1,y,e_1}^*\|}{\|y\|} = \inf_{\theta \in D_2} \|(z(\theta, e_1), \theta_1, \ldots, \theta_{p-1})\|$$
$$> 0 \qquad \text{with probability 1}$$

and, given $X_0 = y$,
$$\lim_{\substack{\|y\| \to \infty \\ y/\|y\| \in D_2}} \|\tilde{\theta}_1 - \tilde{\theta}_1^*\| = \lim_{\substack{\|y\| \to \infty \\ y/\|y\| \in D_2}} \|\theta_{1,y,e_1} - \theta_{1,y,e_1}^*\| = 0 \qquad \text{almost surely}.$$

By the uniform continuity of $\nu$ on $D_2$, therefore,

(5.27)
$$\limsup_{\|x\| \to \infty} E(|\nu(X_{2p+1}) - \nu(\widetilde{X}_{2p+1}^*)| \mathbb{1}_{A_{2p}} | X_0 = x)$$
$$\leq \limsup_{\|x\| \to \infty} E(|\nu(X_{2p+1}) - \nu(\widetilde{X}_{2p+1}^*)| \mathbb{1}_{A_{2p} A_{2p+1}} + 2L_7 \mathbb{1}_{A_{2p+1}^c} | X_0 = x)$$
$$\leq \lim_{\substack{\|y\| \to \infty \\ y/\|y\| \in D_2}} E(|\nu(\tilde{\theta}_1) - \nu(\tilde{\theta}_1^*)| \mathbb{1}_{A_1} | X_0 = y) + 2L_7(2p+1)\varepsilon_3 < \varepsilon_6/4.$$

We may also conclude
$$\limsup_{\substack{\|y\| \to \infty \\ y/\|y\| \in D_2}} E(\log(V_1(X_{1,y,e_1})/V_1(X_{1,y,e_1}^*))) = 0,$$

so that

(5.28) $$\limsup_{\|x\|\to\infty} E(\log(V_1(X_{2p+1})/V_1(\widetilde{X}^*_{2p+1}))\mathbb{1}_{A_{2p}}|X_0=x) = 0.$$

Hence, by (5.26)–(5.28),

(5.29) $$\limsup_{\|x\|\to\infty}|E((\nu(X_{2p+1}) - \nu(X_{2p}) + \log(V_1(X_{2p+1})/V_1(X_{2p})))\mathbb{1}_{A_{2p}}|X_0=x) - \log\rho| < 3\varepsilon_6/4.$$

The definitions of $L_7$ and $K_4$ and (5.25) provide

(5.30) $$\limsup_{\|x\|\to\infty}|E((\nu(X_{2p+1}) - \nu(X_{2p}) + \log(V_1(X_{2p+1})/V_1(X_{2p})))\mathbb{1}_{A^c_{2p}}|X_0=x)| < 2p(2L_7+K_4)\varepsilon_3 \leq \varepsilon_6/4.$$

By (5.29) and (5.30), we conclude

(5.31) $$\limsup_{\|x\|\to\infty}|E(\nu(X_{2p+1}) - \nu(X_{2p}) + \log(V_1(X_{2p+1})/V_1(X_{2p}))|X_0=x) - \log\rho| < \varepsilon_6.$$

Now, finally, we assume $\rho < 1$ and $\varepsilon_6 > 0$ has been chosen so that $\log\rho_1 \stackrel{\text{def}}{=} \log\rho + \varepsilon_6 < 0$. The result follows from Lemma 5.10 with $n = 2p$. $\square$

PROOF OF THEOREM 3.3. Let $\varepsilon_6$, $\nu(x)$, $L_7$ and $V_1(x)$ be as in the proof of Theorem 3.2. Note that $\varepsilon_6$ is arbitrary but $\nu(x)$ and $L_7$ depend on it. On the basis of Lemma 5.9, choose

$$K_5 = \sup_{x\in\mathbb{R}^p} E(|\log(V_1(X_1)/V_1(x))||X_0=x) < \infty.$$

Also, let $n_0 \geq (4pK_5 + 2L_7)/\varepsilon_6$.

Now let $\overline{B}_j = \nu(X_{2p+j}) - \nu(X_{2p+j-1}) + \log(V_1(X_{2p+j})/V_1(X_{2p+j-1}))$. Then, by (5.31),

$$\limsup_{\|x\|\to\infty}|E(\overline{B}_1|X_0=x) - \log\rho| < \varepsilon_6.$$

Since by Lemma 5.9, $\limsup_{\|x\|\to\infty} P(\|X_1\| \leq M|X_0=x) \to 0$ for every $M < \infty$, it is easy to see inductively that, for every $j \geq 1$,

$$\limsup_{\|x\|\to\infty}|E(\overline{B}_{j+1}|X_0=x) - \log\rho|$$
$$\leq \limsup_{\|x\|\to\infty} E(|E(\overline{B}_{j+1}|X_1) - \log\rho||X_0=x)$$
$$\leq \limsup_{\|x\|\to\infty}|E(\overline{B}_j|X_0=x) - \log\rho| < \varepsilon_6.$$

Therefore,

$$\limsup_{\|x\|\to\infty} \frac{1}{n}|E(\nu(X_{2p+n}) - \nu(X_{2p})$$

(5.32)
$$+ \log(V_1(X_{2p+n})/V_1(X_{2p}))|X_0 = x) - \log\rho|$$

$$= \limsup_{\|x\|\to\infty} \left|\frac{1}{n}\sum_{j=1}^n E(\overline{B}_j|X_0 = x) - \log\rho\right| < \varepsilon_6.$$

Note that

$$E(|\log(V_1(X_n)/V_1(x)) - \log(V_1(X_{2p+n})/V_1(X_{2p}))||X_0 = x)$$
$$= E(|\log(V_1(X_{2p+n})/V_1(X_n)) - \log(V_1(X_{2p})/V_1(x))||X_0 = x) \le 4pK_5.$$

Consequently, by (5.32),

$$\limsup_{\|x\|\to\infty} \frac{1}{n}\left|E\left(\log\left(\frac{V_1(X_n)}{V_1(x)}\right)\Big|X_0 = x\right) - \log\rho\right| < \varepsilon_6 + \frac{4pK_5 + 2L_7}{n} < 2\varepsilon_6,$$

for every $n \ge n_0$. Since $\varepsilon_6$ can be arbitrarily small, this proves (3.4) and hence the Lyapounov exponent must be $\log\rho$. □

PROOF OF COROLLARY 3.4. Since a stationary distribution $\Pi$ exists for $\{\theta_t^*\}$, (3.5) implies

$$\log\rho = \int_\Theta \int_\mathbb{R} \log(w(\theta, u)) f(u) \, du \, \Pi(d\theta)$$
$$= \int_\Theta E(\nu(\theta_1^*) - \nu(\theta) + \log w(\theta, e_1)|\theta_0^* = \theta)\Pi(d\theta) < 0.$$

$V$-uniform ergodicity follows by Theorem 3.2. □

PROOF OF THEOREM 3.5. Define $\Theta_{\#\#}$ as before and assume (i) holds. We first verify that $\{X_t\}$ is $V$-uniformly ergodic. We show below that (i) implies (ii) and then that (ii) implies (i) with $\lambda$ continuous on $\Theta_{\#\#}$, bounded and bounded away from 0. Therefore, there is no loss in assuming that $\lambda$ is continuous on $\Theta_{\#\#}$, bounded and bounded away from 0. Let $V_1(x) = 1 + \|x\|^r$. By an argument analogous to that for Theorem 3.2, one may show that

$$\limsup_{\|x\|\to\infty} \frac{E(\lambda(X_{2p+1}/\|X_{2p+1}\|)V_1(X_{2p+1})|X_0 = x)}{E(\lambda(X_{2p}/\|X_{2p}\|)V_1(X_{2p})|X_0 = x)} < 1.$$

The drift condition for geometric ergodicity is therefore satisfied with test function $V(x) = E(\lambda(X_{2p}/\|X_{2p}\|)V_1(X_{2p})|X_0 = x)$. From (5.5) we have

$$E(|a^*(x) + b^*(x)e_1|^r) \ge E(|a^*(x) + b^*(x)e_1|^r \mathbb{1}_{|a^*(x)+b^*(x)e_1|\ge\varepsilon\|x\|})$$
$$\ge \delta_1\|x\|^r,$$

for some $\delta_1 > 0$, and therefore by (1.1) and the assumptions, $E((1+\|X_1\|^r)|X_0 = x) \geq \delta_2(1 + \|x\|^r)$ for some $\delta_2 > 0$. Thus, $E((1 + \|X_{2p}\|^r)|X_0 = x) \geq \delta_2^{2p}(1 + \|x\|^r)$. Since $\lambda(\theta)$ is bounded away from 0, it follows that $V(x) \geq K\|x\|^r$ for some positive constant $K$. Likewise, $V(x) \leq L + M\|x\|^r$ for finite $L, M$.

Next, we show that (i) implies (ii). Let $\beta$ be the left-hand side of (3.6) and choose $L_9 < \infty$ such that $1/L_9 < \lambda(\theta) < L_9$ for all $\theta$. By (3.6) and the Markov property,

$$E\left(\prod_{t=1}^{n}(w(\theta_{t-1}^*, e_t))^r \Big| \theta_0^* = \theta\right)$$

$$\leq L_9 E\left(E(\lambda(\theta_n^*)(w(\theta_{n-1}^*, e_n))^r|\theta_{n-1}^*)\prod_{t=1}^{n-1}(w(\theta_{t-1}^*, e_t))^r \Big| \theta_0^* = \theta\right)$$

$$\leq L_9 \beta E\left(\lambda(\theta_{n-1}^*)\prod_{t=1}^{n-1}(w(\theta_{t-1}^*, e_t))^r \Big| \theta_0^* = \theta\right)$$

$$\leq \cdots \leq L_9 \beta^n \lambda(\theta) \leq L_9^2 \beta^n.$$

The condition in (ii) follows easily.

Finally, we assume (ii) and verify that (i) holds. We need to construct a test function $\lambda(\theta)$ that is bounded, bounded away from 0 and continuous on $\Theta_{\#\#}$. This has two parts.

First, we define $\lambda$ and show it satisfies (3.6). Choose $n \geq 1$ and $\delta > 0$ such that

(5.33) $$\sup_{\theta \in \Theta}\left(E\left(\prod_{t=1}^{n}(\delta + w(\theta_{t-1}^*, e_t))^r \Big| \theta_0^* = \theta\right)\right) < 1$$

and define $Q_t = (\delta + w(\theta_{t-1}^*, e_t))^r$. The test function we need is

$$\lambda(\theta) = \prod_{t=1}^{n-1}(E(Q_t \cdots Q_1|\theta_0^* = \theta))^{1/n}.$$

Applying Hölder's inequality to $n$ factors,

$$E(\lambda(\theta_1^*)(w(\theta, e_1))^r|\theta_0^* = \theta)$$

$$\leq E\left(\prod_{t=1}^{n-1}(E(Q_{t+1}\cdots Q_2|\theta_1^*)Q_1)^{1/n}Q_1^{1/n} \Big| \theta_0^* = \theta\right)$$

$$\leq \prod_{t=1}^{n}(E(Q_t \cdots Q_1|\theta_0^* = \theta))^{1/n} = (E(Q_n \cdots Q_1|\theta_0^* = \theta))^{1/n}\lambda(\theta),$$

thus verifying (3.6) by (5.33).

Second, we show $\lambda$ has the desired properties. Define $q_t(\theta) = E(Q_t \cdots Q_1 | \theta_0^* = \theta)$. By the Markov property,

(5.34) $$q_t(\theta) = E(q_{t-1}(\theta_1^*)Q_1 | \theta_0^* = \theta).$$

Let $K_6 = \sup_\theta E(Q_1 | \theta_0^* = \theta)$ and by iteratively using (5.34) we have $\delta^t \leq q_t(\theta) \leq K_6^t$ and $\delta^{(n-1)/2} \leq \lambda(\theta) \leq K_6^{(n-1)/2}$ for all $\theta$.

Note that $\lambda$ is continuous on $\Theta_{\#\#}$ if each $q_t(\theta)$ is. Clearly, $q_1$ is continuous on $\Theta_{\#\#}$ since $w(\theta, e_1)$ is and $\{(w(\theta, e_1))^r\}_{\theta \in \Theta}$ is uniformly integrable. By the argument of Lemma 5.5, $q_{t-1}(\theta)$ continuous on $\Theta_{\#\#}$ implies $q_t(\theta) = E(q_{t-1}(\theta_1^*)Q_1 | \theta_0^* = \theta)$ is also, and we obtain the conclusion inductively. $\square$

PROOF OF EXAMPLE 4.4. Define $c_1 = b_{11}^r E_1 + b_{21}^r E_2$ and $c_i = (b_{1i}^r p_1 + b_{2i}^r p_2) E(|e_1|^r)$, $i = 2, \ldots, p$. Choose $\beta \in (0,1)$ such that $\sum_{i=1}^p \beta^{-i} c_i = 1$. Now let $d_{jp} = \beta^{-1} b_{jp}^r E(|e_1|^r)$, $j = 1, 2$ and

$$d_{ji} = \beta^{-1} b_{ji}^r E(|e_1|^r) + \sum_{k=i+1}^p \beta^{i-k-1} c_k, \qquad i = 1, \ldots, p-1, \ j = 1, 2.$$

Thus, $d_{11} E_1 + d_{21} E_2 = E(|e_1|^r)$,

(5.35) $$b_{jp}^r(d_{11} E_1 + d_{21} E_2) = \beta d_{jp}, \qquad j = 1, 2,$$

and for $i = 2, \ldots, p-1$, $j = 1, 2$,

(5.36) $$b_{ji}^r(d_{11} E_1 + d_{21} E_2) + (d_{1,i+1} p_1 + d_{2,i+1} p_2) E(|e_1|^r) = \beta d_{ji}.$$

Now define $d_{j,p+1} = 0$, $j = 1, 2$, $R_1 = \{x : x_1 \leq 0\}$, $R_2 = \{x : x_1 > 0\}$ and

$$V(x) = \begin{cases} d_{11}|x_1|^r + \cdots + d_{1p}|x_p|^r, & \text{if } x_1 \leq 0, \\ d_{21}|x_1|^r + \cdots + d_{2p}|x_p|^r, & \text{if } x_1 > 0. \end{cases}$$

Applying (5.35) and (5.36),

$$E(V(X_1^*) | X_0^* = x)$$

$$= \sum_{j=1}^2 \left( \sum_{i=1}^p b_{ji}^2 |x_i|^2 \right)^{r/2} \mathbb{1}_{R_j}(x)(d_{11} E_1 + d_{21} E_2)$$

$$+ \sum_{j=1}^2 \sum_{i=2}^p (d_{1,i+1} p_1 + d_{2,i+1} p_2) |x_i|^r \mathbb{1}_{R_j}(x) E(|e_1|^r)$$

$$\leq \sum_{j=1}^2 \sum_{i=1}^p (b_{ji}^r(d_{11} E_1 + d_{21} E_2) + (d_{1,i+1} p_1 + d_{2,i+1} p_2) E(|e_1|^r)) |x_i|^r \mathbb{1}_{R_j}(x)$$

$$= \sum_{j=1}^2 \sum_{i=1}^p \beta d_{ji} |x_i|^r \mathbb{1}_{R_j}(x) = \beta V(x),$$

with equality when $r = 2$.

Letting $\lambda(\theta) = V(\theta)$, for $\theta \in \Theta$, leads immediately to (3.6) and the conclusion.

□

DEPARTMENT OF STATISTICS
TEXAS A&M UNIVERSITY
COLLEGE STATION, TEXAS 77843-3143
USA
E-MAIL: dcline@stat.tamu.edu
E-MAIL: pu@stat.tamu.edu